\documentclass[11pt]{article}
\usepackage{amsfonts}
\usepackage{enumerate}
\usepackage{makeidx}
\usepackage{amsthm,amssymb}
\usepackage{latexsym}
\usepackage{graphicx}
\usepackage{amsmath}
\usepackage[T1]{fontenc}
\usepackage[utf8]{inputenc}
\usepackage{amsfonts}
\usepackage[normalem]{ulem}

\newtheorem{thm}{Theorem}[section]
\newtheorem{lem}{Lemma}[section]
\numberwithin{equation}{section}
\theoremstyle{definition}

\newtheorem{rem}{Remark}[section]

\newcommand{\N}{\mathbb{N}}

\newcommand{\E}{\mathbf{E}}
\newcommand{\V}{\mathbf{Var}}
\newcommand{\R}{\mathbb{R}}
\newcommand{\C}{\mathbb{C}}
\def\Cov{\mathbf{Cov}}
\newcommand{\al}{\alpha}

\newcommand{\Y}{Y}

\newcommand{\YY}{\mathcal{Y}_n}
\newcommand{\y}{\mathbf{y}}

\def\Tr{{\mathrm{Tr}}}

\makeatother

\begin{document}
\title{A note on the fluctuations of the resolvent traces of a tensor model of sample covariance matrices}
\author{Alicja Dembczak-Kołodziejczyk}
\newcommand\address{\noindent\leavevmode

\noindent
Alicja Dembczak-Kołodziejczyk,\\
University of Opole,\\
48 Oleska, 45-052,\\
Opole, Poland.\\
\texttt{\small
e-mail:
alicja.dembczak@uni.opole.pl
}
}
\date{}
\maketitle

\begin{abstract}
\noindent

In this note, we consider a sample covariance matrix of the form  $$M_{n}=\sum_{\al=1}^m \tau_\alpha {\y}_{\alpha}^{(1)} \otimes {\y}_{\alpha}^{(2)}({\y}_{\alpha}^{(1)} \otimes {\y}_{\alpha}^{(2)})^T,$$ where $(\y_{\alpha}^{(1)},\, {\y}_{\alpha}^{(2)})_{\al}$ are independent vectors uniformly distributed on the unit sphere $S^{n-1}$ %copies of a random vector $\y\in S^{n-1}$ %be a random vector uniformly distributed on the unit sphere in
%\R^n$, whose distribution satisfies $\E y_i=0, \E y_iy_j=n^{-1}\delta_{ij}$, where $i,j\le n$
 and $\tau_\alpha \in \R_+ $.
%Given a smooth enough test function $\varphi$ ($\varphi \in \mathcal{H}_s$ for all $s>6.5$),
We show that as $m, n \to \infty$, $m/n^2\to c>0$, the centralized traces of the resolvents, $\Tr(M_n-zI_n)^{-1}-\E\Tr(M_n-zI_n)^{-1}$, $\Im z\ge \eta_0>0$, converge in distribution to a two-dimensional
Gaussian random variable with zero mean and a certain covariance matrix. This work is a continuation of \cite{DKL:23, Ly:18}. % the linear eigenvalue statistics $\Tr \ \varphi(M_n)-\E \Tr\varphi(M_n)$.
%We %prove the CLT for linear eigenvalue statistic and  show for
%how that given a test function $\varphi$ ($\varphi \in \mathbf{H}^s$ for all $s>5.5$), the  linear eigenvalue statistic $\Tr \varphi(M_n)$  converges in distribution to a gaussian random variable.\\% with zero mean and variance which we denote by  $V[\varphi]$.
\\


\end{abstract}
\section{Introduction: Model and Main Results}

 Given  $n \in\N$ and $m=m(n)\in \N$, consider  an $n^2\times n^2$ sample covariance matrix of the form
 \begin{equation}\label{M}
   M_{n}=\sum_{\al=1}^m \tau_\alpha Y_\al Y_\al^T,
 \end{equation}
where for every $n\in N, \,\{\tau_\al: \al=1,\ldots, m\}$ are nonnegative numbers,
\begin{equation}\label{Y}
Y_\al=\y_\al^{(1)}\otimes \y_\al^{(2)}=(\y_{\al 1}^{(1)}\y_{\al 1}^{(2)},\dots,\y_{\al n}^{(1)}\y_{\al n}^{(2)} )\in S^{n^2-1},
\end{equation}
and $(\y_{\alpha}^{(1)},\, {\y}_{\alpha}^{(2)})_{\al=1}^m$ are independent vectors uniformly distributed on the unit sphere $S^{n-1}$.
Denote $(\lambda_j^{(n)})_{j=1}^{n^2}$ the eigenvalues of $M_n$ counting their multiplicity, and introduce their normalized counting measure  $N_n$ by setting for any $\Delta\subset \R$
\begin{align*}
N_n(\Delta)=|\{\lambda_j\in \Delta, j = 1, . . . , n\}|/n^2.
\end{align*}
Similarly, we define the normalized counting measure $\sigma_n=\sigma_{m(n)}$ of $(\tau_\alpha)_{\alpha=1}^{m}$, $\sigma_m(\Delta)=|\{\tau_\alpha\in \Delta, \alpha = 1,\ldots,m\}|/m $,
and assume that the sequence $(\sigma_n)_{n=1}^\infty$ converges weakly to a probability measure $\sigma$.
%Similarly, we define the NCM $\sigma_m$ of $(\tau_\alpha)_{\alpha=1}^m$, $\sigma_m(\Delta)=|\{\tau_\alpha\in \Delta, \alpha \in[1,m]\}|/m $,
%and assume that the sequence $(\sigma_m)_{m=1}^\infty$ converges weakly to a probability measure $\sigma$,
%\begin{align}\label{tau}
%\lim\limits_{m\to\infty}\sigma_m=\sigma, \quad \text{d}|\sigma_m-\sigma|=O(n^{-2}), \,\,\, \sigma(\R)=1.
%\end{align}
It is known (see \cite{Ly:18, Pa-Pa:09}) that as $m, n \to \infty$, $m/n^2\to c >0$, the measure $N_n$ converges weakly with probability 1 to the probability measure
$N$ which Stieltjes transform
\begin{align*}
f(z):=\int_\R\frac{N(d\lambda)}{\lambda-z},\quad \Im z \neq 0,
\end{align*}
 is uniquely determined by the functional equation
\begin{align}\label{rownanie}
zf(z)=-1+c-\int\frac{c\ \text{d}\sigma(\tau)}{1+\tau f(z)},
\end{align}
considered in the class of analytic functions in $\C\backslash\R$  such that $\Im f(z)\Im z \geq 0, \Im z\neq 0$. In
particular, introducing the linear eigenvalue statistic $\Tr \, \varphi(M_n)$
%\begin{align*}
%\mathcal{N}_n[\varphi]=\sum_{j=1}^n\varphi(\lambda_j^{(n)})
%\end{align*}
corresponding to a continuous bounded test function $\varphi$,
we get with probability 1
\begin{align*}
\lim\limits_{n\to \infty}n^{-2}\Tr\ \varphi(M_n)=\int\varphi(\lambda)N(d\lambda).
\end{align*}
Model (\ref{M}) -- (\ref{Y}) is a particular case of the model considered in \cite{Ly:18} (we denote it here by $M_n^{(1)}$), where it was assumed that $(\y_\al^{(1)},\y_\al^{(2)})_{\al=1}^m$ were iid copies of a vector $\y=(y_1,\ldots,y_n)$ with an unconditional distribution satisfying  the following moment conditions
%as $n\to \infty$:
\begin{align*}
&\E y_i=0,\,\,\,\E y_iy_j=n^{-1}\delta_{ij}, \,\,\, i,\ j \le n,\\
&a_{2,2}:=\E y_i^2 y_j^2=n^{-2}+an^{-3}+O(n^{-4}), \,\,\,\forall i\neq j,\\
&\E y_j^4-3a_{2,2}=bn^{-2}+O(n^{-3}), \quad n\to \infty
\end{align*}
for some $a, \ b \in \R$.
In \cite{Ly:18}, a central limit theorem (CLT) for the linear eigenvalue statistics of  $M_n^{(1)}$ was proven. It was shown that  for any test function $$\varphi \in\mathcal{H}_s=\Big\{\varphi:\int|\hat{\varphi}(t)|^2(1+t)^{2s}\text{d}t<\infty\Big\},\quad s>5/2,$$  the variance of $\Tr\,\varphi(M_n^{(1)})$ grows to infinity not faster than $n$,
%\begin{equation*}
%  \V \,\Tr \varphi(M_n)=O(n),
%\end{equation*}
and $(\Tr\,\varphi(M_n^{(1)})-\E \Tr\,\varphi(M_n^{(1)}))/\sqrt{n}$ converges in distribution  to  a  Gaussian random variable with zero mean and variance \begin{align*}
V^{(1)}[\varphi]=\lim\limits_{\eta \downarrow 0}\frac{2(a+b+2)c}{\pi ^{2}}\int \tau^2\left(\Im \int\frac{f^\prime(\lambda+i\eta)}{(1+\tau f(\lambda+i\eta))^2}\varphi(\lambda)\text{d}\lambda\right)^2\text{d}\sigma(\tau) %=\lim_{n\to \infty} \V \,\Tr \varphi(M_n)/\sqrt{n}
\end{align*}%
(see Theorem 1.9 of \cite{Ly:18} for the details.)
Note that $a+b+2=0$ for $\y_\al^{(1)},\y_\al^{(2)} \sim U(S^{n-1})$, so that to get a nontrivial limit one needs to renormalise the linear statistic. This was the starting point in \cite{DKL:23}, where the authors dealing with matrices of the form (\ref{M}) -- (\ref{Y}) with $\tau_1=...=\tau_m=1$, considered traces of the resolvent $\gamma_n (z) = \Tr\,(M_n-zI_n)^{-1}$, and shown that in the case of vectors uniformly distributed on the unit sphere the variance of  $\gamma_n$ is not of order n but of order 1. Also  there was found the limit of the covariance of the resolvent traces. Here we continue this research and prove a CLT for the resolvent traces $ \gamma_n  $. %of $M_n$
Before stating the main result, we impose additional conditions on $\tau_\al$.\\
\textbf{Condition A}. (i) We suppose that  $\{\tau_\al\}_{\al}$ are uniformly bounded in $\al, n$:
$$\exists T>0\quad \sup\limits_{\al, n} \tau_\al\le T,$$
(ii) the counting measure $\sigma_n$ of $\{\tau_\al\}_{\al}$ converges weakly to a probability measure $\sigma$, and there exists an absolute constant $C>0$ such that for every bounded continuous function $\varphi$
$$\Big|\int\varphi(\tau)\text{d}\sigma_n-\int\varphi(\tau)\text{d}\sigma
\Big|\le\frac{C}{n^2}\sup\limits_\tau|\varphi(\tau)|.$$

Our main result is the following theorem.
\begin{thm}\label{tw} Given $m,n\in \N$, % Let  $m/n^2\to c$, $Y$ be a vector satysfying (\ref{a22}).
let $M_n$ be defined in (\ref{M}) -- (\ref{Y}), where $\{\tau_\al\}_{\al}$ satisfy Condition A.
Let $\gamma_n(z)=\Tr \, (M_n-zI_n)^{-1}$, $z=\xi +i\eta \in \C\backslash\R$.
Then there exists $\eta_0>0$ such that for every
$
z \in \C_{\eta_0}:=\{z\in \C:\eta>\eta_0\}
$
the centered trace of the resolvent, $\gamma_n-\E\gamma_n$,  converges in distribution as  $m/n^2\to c>0$, $m,n\to \infty$, to a  two-dimensional Gaussian random variable with zero mean and a covariance matrix of the form:
\begin{align}\label{mkow}
\Sigma=\begin{pmatrix}\frac{1}{4}(2C(z,\overline{z})+C(z,z)+C(\overline{z},\overline{z}))& \frac{1}{4i}(C(z,z)- C(\overline{z},\overline{z}))\\ \frac{1}{4i}(C(z,z)- C(\overline{z},\overline{z}))&\frac{1}{4}(2C(z,\overline{z})-C(z,z)-C(\overline{z},\overline{z}))\end{pmatrix},
\end{align}
where
%and
\begin{align}
%&V(z)%[\varphi]
%=\frac{1}{4}[2C(z,\overline{z})-C(z,z)-C(\overline{z},\overline{z})],\\
C(z_1,z_2)%:=\lim\limits_{n\to\infty}\Cov\{\gamma(z_1),\gamma(z_2)\}
=\frac{\partial ^{2}}{\partial z_1\partial z_2}\bigg(
2\log \frac{\Delta f}{\Delta z}+2f(z_1)f(z_2)\frac{\Delta z}{\Delta f}+3\int\frac{c^2\tau^4f^2(z_1)f^2(z_2)\text{d}\sigma(\tau)}{(1+\tau f(z_1))^2(1+\tau f(z_2))^2}\bigg)\label{CZZ},%\Big(\frac{cf^2(z_1)}{(1+f(z_1))^2}\Big)'\Big(\frac{cf^2(z_2)}{(1+f(z_2))^2}\Big)'\bigg)\label{CZZ}
\end{align}
$\Delta f=f(z_1)-f(z_2)$, $\Delta z =z_1-z_2$, %$z_1, z_2\in\{z\in \C :\Im z\geq \eta_0\}$
and $f$ is defined in  (\ref{rownanie}).%\cite{Ma-Pa:67}.
\end{thm}

%\begin{lem}\label{lem} Under the  Conditions of Theorem \ref{tw}. %%let $M_n$ be defined in (\ref{M})--(\ref{Y}), where $\{\tau_\al\}_{\al}$ satisfy Condition A.
%Let $\gamma_n(z)=\Tr \, (M_n-zI_n)^{-1}$, $ z\in \C\backslash\R$.
%Then there exists $\eta_0>0$ such that for every
%\begin{align*}
%%z \in \C_{\eta_0}:=\{z\in \C:\Im z>\eta_0\}
%\end{align*}
%in the limit $m=m(n)\to \infty$ and  $m/n^2\to c>0$ as $n\to \infty$
%We have $$\Im (\gamma_n(z) -\E\gamma_n(z))\xrightarrow[n\to \infty]{d}Z\sim N(0,V(z)),$$ % , $\mathcal{N}_n^\circ[\varphi]$
%converges in distribution to the Gaussian random variable with zero mean and variance $V%[\varphi]
%=\lim_{\eta\downarrow 0} V_\eta%[\varphi]
%$,
%where
%\begin{align}
%&V(z)%[\varphi]
%=\frac{1}{4}[2C(z,\overline{z})-C(z,z)-C(\overline{z},\overline{z})].
%&C(z_1,z_2)%:=\lim\limits_{n\to\infty}\Cov\{\gamma(z_1),\gamma(z_2)\}
%=\frac{\partial ^{2}}{\partial z_1\partial z_2}\bigg(
%2\log \frac{\Delta f}{\Delta z}+2f(z_1)f(z_2)\frac{\Delta z}{\Delta f}+3\int\frac{c^2\tau^4f^2(z_1)f^2(z_2)\text{d}%\sigma(\tau)}{(1+\tau f(z_1))^2(1+\tau f(z_2))^2}\bigg)\label{CZZ2},%\Big(\frac{cf^2(z_1)}{(1+f(z_1))^2}\Big)'\Big(\frac{cf^2(z_2)}{(1+f(z_2))^2}\Big)'\bigg)\label{CZZ}
%\end{align}
 %$\Delta f=f(z_1)-f(z_2)$, $\Delta z =z_1-z_2$, %$z_1, z_2\in\{z\in \C :\Im z\geq \eta_0\}$
% and $f$ is defined in  (\ref{rownanie}).%\cite{Ma-Pa:67}.
%\end{lem}
\begin{rem}
(i) In (\ref{CZZ}),
$
C(z_1,z_2)=\lim\limits_{n\to\infty}\Cov\{\gamma_n(z_1),\gamma_n(z_2)\}.%=\lim\limits_{n\to\infty}\E\{\gamma_n(z_1)\gamma_n^\circ(z_2))\},
$

%where $\gamma_n(z)=\Tr (M_n-zI_n)^{-1}$ and $\gamma_n^\circ=\gamma_n-\E \gamma_n$.\\
%(ii) It can be shown that if $\tau_1=\ldots=\tau_m=1$ then (\ref{t:clt})--(\ref{CZZ}) yield (\ref{limitvariance}).\\
(ii) Condition A(i) is a technical one, as (\ref{CZZ}) points out, it should be  enough for $\sigma$ to have two finite moments.
\end{rem}
%For $\tau_\alpha = 1$ we get our result coinsides with the result of Theorem 1.1 in \cite{DKL:23}.
%\begin{align*}
%V[\varphi ]&=\frac{1}{2\pi ^{2}}\int_{a_{-}}^{a_{+}}%
%\int_{a_{-}}^{a_{+}} \left(\frac{\triangle \varphi }{\triangle \lambda }%
%\right) ^{2}\frac{(4c-(\lambda _{1}-a_{m})(\lambda
%_{2}-a_{m}))d\lambda
%_{1}d\lambda _{2}}{\sqrt{(a_{+}-\lambda _{1})(\lambda _{1}-a_{-})}\sqrt{%
%(a_{+}-\lambda _{2})(\lambda _{2}-a_{-})}}
%\\
%&-\frac{2}{4c\pi ^{2}}\left( \int_{a_{-}}^{a_{+}}\varphi (\mu )\frac{\mu
%-a_{m}}{\sqrt{(a_{+}-\mu )(\mu -a_{-})}}d\mu \right) ^{2}+\frac{3}{(c\pi) ^{2}}\left( \int_{a_{-}}^{a_{+}}\varphi (\mu )\frac{2c-(\mu-a_m)^2}{\sqrt{(a_{+}-\mu )(\mu -a_{-})}}d\mu \right) ^{2}.\notag
%\end{align*}
%where
%$$a_{\pm }=(1\pm \sqrt{c})^{2},\quad a_{m}=(a_++a_-)/2=1+c$$
%(see Theorem 1.1  and formulas (1.6), (1.8) of \cite{DKL:23} for the details.)\\

%The proof of the CLT for linear eigenvalue statistics (proof of Theorem \ref{tw}) based on the scheme proposed in \cite{Sh:11}.

%Our research is based on the proof of CLT, and show  for  a  smooth test function $\varphi$ ($\varphi \in \mathbf{H}^s$ for all $s>5.5$), the  linear eigenvalue statistic  $\Tr \varphi(M_n)$ converges in distribution to a gaussian random variable. In Section 2  there are notations. In Section 3 is the variance of the linear eigenvalue statistic, which is needed to prove CLT. Section 4 contains the proof of Theorem \ref{tw}.

This result is an analog of Bai and Silverstein's master lemma {\cite[Lemma~1.1]{Ba-Si:04}},  it can be used to prove the CLT for an arbitrary smooth enough test function. In Section 2 we introduce the main notations and present some auxilary results, then in Section 3 we give a proof of the theorem. Our proof of Theorem \ref{tw} heavily relies on \cite{DKL:23, Ly:18} and follows the scheme proposed in \cite{Pa-Sh, Sh:11}.

\section{Notations and preliminary results}
Let ${\y}^{(1)}$ and ${\y}^{(2)}$  be  i.i.d. copies of $\y\sim U(S^{n-1})$, and let $\Y={\y}^{(1)} \otimes {\y}^{(2)}=\big(y_{i}^{(1)} y_{j}^{(2)}\big)_{i,j=1}^n.$  In what follows, without loss of generality, we will assume that $m=cn^2$.
It is easy to show  that
\begin{align}
&\E Y_{ij}=0,\quad \E Y_{ij}Y_{pq}=n^{-2}\delta_{ij}\delta_{pq},\notag
\\
&\E Y_{ij}^2Y^2_{pq}=a^2_{2,2}=\frac{1}{n^2(n+2)^2}=\frac{1}{n^4}-\frac{4}{n^5}+\frac{12}{n^6}+O(n^{-7}), \label{a22}
\\
&\E Y_{jj'}Y_{ss'}Y_{pp'}Y_{qq'}= a^2_{2,2}(\delta_{js}\delta_{pq}+\delta_{jp}\delta_{sq}+\delta_{jq} \delta_{sp})
(\delta_{j's'}\delta_{p'q'}+\delta_{j'p'}\delta_{s'q'}+\delta_{j'q'} \delta_{s'p'})\notag
\end{align}
(see \cite{DKL:23}.) Note that if $H$ is an arbitrary matrix and $Y$ is defined in (\ref{Y}), then by (\ref{a22})
\begin{align}\label{wariancjah}
	\V(HY,Y)&=(a_{2,2}^2-n^{-4})|\Tr\ H|^2+2a_{2,2}^2\Tr\ |H|^2+4a_{2,2}^2\sum_{j,s,p,q=1}^nH_{js,ps}\overline{H}_{jq,pq}\\&\quad+2a_{2,2}^2\sum_{j,s,p,q=1}^nH_{js,pq}\overline{H}_{ps,jq},\notag
\end{align}
where $\Tr |H|^2=\Tr HH^*$  and without loss of generality we assume that
$$\sum_{j,s,p,q} H_{js,ps}\overline{H}_{jq,pq}=\sum_{j,s,p,q} H_{sj,sp}\overline{H}_{qj,qp}$$
(they are asymptotically equal and also can be treated in the same way.) By the Cauchy-Schwarz inequality $|\Tr\, H|^2\le n\Tr |H|^2$, $|\sum H_{js,ps}\overline{H}_{jq,pq}|\le n \Tr |H|^2$, and $|\sum H_{js,pq}\overline{H}_{ps,jq}|\le \Tr |H|^2$, so that
\begin{align}\label{VH}
	\V(HY,Y)\lesssim n^{-1}\|H\|_{op},
\end{align}
where $\|H\|_{op}$ denotes the operator norm of $H$. Here and in what follows  we use notation $|A|\lesssim |B|$ if there exists an absolute constant $C$ or a positive quantity, which does not depend on the parameters of the model, such as $n$,$m$, $\tau_\al$, $Y_\al$, but possibly depends on $\eta_0$, $T$, or the order of the corresponding moment.

Given $z \in \C\backslash\R$, let $G(z)=(M_n-zI_n)^{-1}$ be the resolvent of $M_n$, and
\begin{align}
&\gamma_n:=\Tr \ G,\quad g_n:=n^{-2}\gamma_n,\quad f_n:=\E g_n.\label{gn}
\end{align}
Our proof  is heavily based on the asymptotic analysis of (\ref{wariancjah}), where $H=G(z)$ is the resolvent of $M_n$ independent of $Y$. In particular, we will need to understand the asymptotic behavior of $g_n$ (and $\gamma_n$, $f_n$) and also of the resolvent counterparts of the two last terms on the r.h.s of (\ref{wariancjah}),
\begin{align}
%\gamma_n(z)=\Tr\ G(z), \,\,\,
&g_n^{(1)}(z_1,z_2):=\frac{1}{n^3}\sum_{j,s,p,q}G_{js, ps}(z_1)G_{jq,pq}(z_2), \quad (f_n^{(1)}:=\E g_{n}^{(1)}),\quad  \text{and} \notag\\
&g_n^{(2)}(z_1,z_2):=\frac{1}{n^2}\sum_{j,s,p,q}G_{js, pq}(z_1)G_{ps,jq}(z_2), \quad (f_n^{(2)}:=\E g_{n}^{(2)}),\label{fn2}
\end{align}
(see Lemma \ref{lematwariancja} below.)\\

 Let $M_{n}^\al=M_{n}- \tau_\alpha Y_\al Y_\al^T=\sum_{\beta\neq\al} \tau_\beta Y_\beta Y_\beta^T$.
In what follows, we use the upper index $\al$ for the quantities corresponding to $M_n^\al$:
$$ G^\al(z):=(M^\al_{n}-z)^{-1},\quad \gamma_n^\al:=\Tr \ G^\al,\quad g_n^\alpha: =n^{-2}\Tr\ G^\alpha,\quad f_n^\alpha:=\E g_n^\alpha.$$
%Introduce $n\times n$ matrix $\G$ by the formula
%\begin{align*}
%\G=\Big(\sum_s G_{js,ps}\Big)_{j,p}
%\end{align*}
Note that $G$ and $G^\alpha$ are connected by the so-called rank one perturbation formula
 \begin{equation}
 G-G^\alpha=-\frac{\tau_\alpha G^\alpha \Y_{\alpha}\Y_\al^T G^\alpha}{1+\tau_\alpha(G^\alpha Y_{\alpha },Y_{\alpha })},\label{G-G}
\end{equation}
which allows to separate $\Y_\alpha$ from the rest of the vectors. Let  $\E_\al=\E_{Y_\al}$ denote the expectation with respect to $Y_\al$. By (\ref{a22}),
\begin{align}\label{eega}
\E_\al(G^\al Y_\al,Y_\al)=g_n^\al,\,\,\,\E(G^\al Y_\al,Y_\al)=f_n^\al.
\end{align}
%g_n^{(1)\al}(z_1,z_2)=n^{-3}\Tr \,\G^\al(z_1)\G^\al(z_2),
 %Let  $\E_\al=\E_{Y_\al}$ denote expectation with respect to $Y_\al$. We have by (\ref{a22}),
%$$
%\E_\al(G^\al Y_\al,Y_\al)=g_n^\al,\,\,\,\E(G^\al Y_\al,Y_\al)=f_n^\al.
%$$
Also given $\xi=\xi(Y_1,\dots,\Y_m)$, we put
$$
\xi^\circ=\xi-\E\xi,\quad
(\xi)_\al^\circ=\xi-\E_\al\xi,
$$
so that $\V \,\xi=\E|\xi^\circ|^2$, $\V_\al \,\xi=\E_\al|\xi-\E_\al\xi|^2=\E_\al|(\xi)_\al^\circ|^2$.
Note that since
\begin{align}\label{rozpisksi}
\frac{1}{\xi}=\frac{1}{\E \xi}-\frac{\xi^\circ}{\xi\E\xi}
\end{align}
for every $k\in \N$ we have
\begin{align}
\frac{1}{\xi}=\frac{1}{\E\xi}-\frac{\xi^\circ}{(\E\xi)^2}+\ldots+\frac{(-1)^{k-1}\xi^{\circ^{k-1}}}{(\E\xi)^{k}}+\frac{(-1)^k\xi^{\circ^k}}{\xi(\E\xi)^{k}}.\label{rozpis}
\end{align}
If $|\xi|$ and $|\E \xi|$ are bounded from below by a non-zero constant, then as  it follows from (\ref{rozpisksi}), 
\begin{align*}
\V \frac{1}{\xi}%=\V\left(\frac{1}{\E\xi}-\frac{\xi^\circ}{\xi \E\xi}\right)
=\frac{1}{(\E\xi)^2}\V \frac{\xi^\circ}{\xi}%=\frac{1}{(\E\xi)^2}\E\Big|\Big(\frac{\xi^\circ}{\xi}\Big)^\circ\Big|^2
\lesssim \V \xi ,
\end{align*}
 Also we have similar elementary inequalities for the central moments of the products and quotients of two bounded random variables. Namely, if $\zeta$, $\eta$ are bounded and independent then $
\V \zeta \eta=\E\zeta^2\E\eta^2-(\E\zeta)^2(\E\eta)^2=\E\zeta^2\V\eta+(\E\eta)^2\V\eta\lesssim \V \zeta +\V \eta$
and more generally
\begin{align*}
\E|(\zeta \eta)^\circ|^p\lesssim \E|\zeta^\circ\eta^\circ|^p+\E|\zeta^\circ\E\eta|^p+\E|\eta^\circ\E\zeta|^p\lesssim \E|\zeta^\circ|^p+\E|\eta^\circ|^p,
\end{align*}
which implies
\begin{align}\label{gwiazdka2.4}
\E\Big|\Big(\frac{\eta}{\zeta}\Big)^\circ\Big|^p\le\E\Big|\frac{\eta}{\zeta}-\frac{\E \eta}{\E\zeta}\Big|^p=
 \E\Big|\frac{\eta^\circ}{\E\zeta}-\frac{ \eta}{\zeta}\cdot \frac{ \zeta^\circ}{\E\zeta}\Big|^p\lesssim \E|\zeta^\circ|^p+\E|\eta^\circ|^p.
\end{align}
In Lemma \ref{lematwariancja} below, we summarize some properties of $\gamma_n$, $f_n^{(1,2)}=\E g_n^{(1,2)}$, and some related quantities that we will need in the following. Before formulating the lemma, we will derive a few more auxiliary inequalities and estimates.
%In {\cite[Lemma~1.1]{DKL:23}}, theorems analogous to those given below were proven for the case when $\tau_1=\ldots=\tau_m=1$ and can be easily extended to the general case (\ref{M}) -- (\ref{Y}) with the appropriate modifications. Here, we will only outline some important steps. 
Note first that, as it follows from (\ref{G-G})
\begin{align}\label{gn-gna}
		\gamma_n-\gamma_n^\alpha=-\frac{\tau_\al(G^{\alpha 2}Y_{\alpha},Y_{\alpha})}{1+\tau_\al(G^\alpha Y_{\alpha },Y_{\alpha })}=:-\frac{B_\al}{A_{\alpha}}.%\quad \text{and}
	\end{align}
It was shown in \cite{DKL:23} (see also \cite{Sh:11}) that, under the assumption $\|Y_\al\|_2=1$, and Condition A(i)% and $|f(z)|\le |\Im z|^{-1}$
, $A_\al$ satisfies
\begin{align}\label{|EA|}
1\lesssim |A_\al|,\ |\E_\al A_\al|, \ |\E A_\al|\lesssim 1,
\end{align}
where the inequalities hold
uniformly in $z\in \C_{\eta_0}$, and so does $B_\al$.
According to the spectral theorem for real symmetric matrices, there exists a non-negative measure \( m_\alpha \) such that
\begin{align*}
(G_\alpha Y_\alpha, Y_\alpha) = \int  \frac{m_\alpha(d\lambda)}{\lambda - z}\quad \text{and}\quad
(G_\alpha^2 Y_\alpha, Y_\alpha) = \int \frac{m_\alpha(d\lambda) }{(\lambda - z)^2}.
\end{align*}
This yields
\begin{align*}
|A_\alpha| \geq |\Im A_\alpha| = |\tau_\alpha|\eta\int\frac{m_\alpha(d\lambda) }{|\lambda - z|^2}\ge \eta |B_\al|,
\end{align*}
implying that $|B_\al/A_\al|\le \eta^{-1}$ and
\begin{align}\label{B/A}
|\gamma_n-\gamma_n^\al|\le \eta_0^{-1}
\end{align}
for any $z\in \C_{\eta_0}$. Also, it follows from (\ref{VH}) and $\|G(z)\|_{op}\le\eta^{-1}$, that
\begin{align}\label{wariancjaab}
\V_\al A_\al, \V_\al B_\al=O(n^{-1}).
\end{align}
This allows us to show that
\begin{align}\label{Vgamman}
		\V \gamma_n(z)=O(n), \quad n\to \infty.
	\end{align}
Indeed, applying (\ref{gwiazdka2.4}) -- (\ref{|EA|}) and (\ref{wariancjaab}) we get
\begin{align*}
		\V\gamma_n&\le \sum_{\alpha}\mathbf{E}|\gamma_n-\E_\al\gamma_n|^{2}
=\sum_{\al}\E|\gamma_n-\gamma_n^\al-\E_\al(\gamma_n-\gamma_n^\al)|^2
	\notag\\&\le\sum_\al\E\Big|\frac{B_\al}{A_\al}-\frac{\E_\al  B_\al}{\E_\al  A_\al}\Big|^2
= \sum_\al  \E\Big|\frac{B_\al^\circ}{\E_\al A_\al}-\frac{B_\al}{A_\al}\cdot \frac{A_\al^\circ}{\E_\al A_\al}\Big|^2\\&
\lesssim \sum_\al \E (\V_\al A_\al+\V_\al B_\al)=O(n)
\end{align*}
(see also Proposition 2 in \cite{Sh:11} and Lemma 3.2 in \cite{GLPP:13}.) This together with (\ref{wariancjaab}) -- (\ref{Vgamman}), the equality $A_\al^\circ=(A_\al)_{\al}^\circ +\tau_\al g_n^{\al \circ}$, and Condition A(i) yields
\begin{align}\label{VAB}
\V A_\al\lesssim\E \V_\al A_\al +\V g_n^\al=O(n^{-1})
\end{align}
and
$\V B_\al =O(n^{-1})$.
More careful and tedious calculations based on the obtained results and analysis of $f_n, f_n^{(1)}, f_n^{(2)}$ of (\ref{gn}) -- (\ref{fn2}) allow us to improve these bounds and show that $\V A_\al=O(n^{-2})$, $\V \gamma_n=O(1)$. To this end we extend Lemma 3.1 in \cite{DKL:23} to our model (\ref{M}) -- (\ref{Y}), and we have:

\begin{lem}\label{lematwariancja}
Under the conditions of Theorem \ref{tw},
\begin{align*}
&(i)\,\,\, \V A_\al, \ \V B_\al=O(n^{-2}),\quad \V \gamma_n=O(1),\\
&(ii)\,\,\,
f_n(z)=f(z)+O(n^{-2}),\\
&(iii)\,\,\, f_n^{(1)}(z_1,z_2)=f(z_1)f(z_2)+\frac{1}{n}\int\frac{c\tau^2 f^2(z_1)f^2(z_2)\text{d}\sigma(\tau)}{(1+\tau f(z_1))(1+\tau f(z_2))}+O(n^{-3/2}),\\
&(iv) \,\,\, f_n^{(2)}(z_1,z_2)=f(z_1)f(z_2)+\int\frac{c\tau^2 f^2(z_1)f^2(z_2)\text{d}\sigma(\tau)}{(1+\tau f(z_1))(1+\tau f(z_2))}+O(n^{-1/2}).
\end{align*}
\end{lem}
\begin{proof}
The proof of Lemma \ref{lematwariancja} repeats with corresponding changes the proof of Lemma 3.1 in \cite{DKL:23}.  Note first that
\begin{align*}
\V A_\al =\tau_\al^2K_n(z,\overline{z}),
\end{align*} where by (\ref{eega}) and (\ref{Vgamman})
\begin{align}\label{Dn}
K_n(z_1,z_2)&:=\Cov\{(G^\al(z_1) Y_\al,Y_\al),(G^\al(z_2) Y_\al,Y_\al)\}\notag\\&=\E(G^\al(z_1) Y_\al,Y_\al)_\al^\circ(G^\al(z_2) Y_\al,Y_\al)+\Cov\{ g_n^\al(z_1),g_n^\al(z_2)\}\notag\\&=\E (G^\al(z_1) Y_\al,Y_\al)_\al^\circ(G^\al(z_2) Y_\al,Y_\al)+O(n^{-3}).
\end{align}
Applying (\ref{a22}) -- (\ref{wariancjah}) we get
\begin{align}\label{nCovAa}
 \E\E_\al (G^\al(z_1) Y_\al,Y_\al)_\al^\circ&(G^\al(z_2) Y_\al,Y_\al)=\frac{1}{n}\E\Big[\frac{-4n^2-4n}{(n+2)^2} g_n^\al(z_1)g_n^\al(z_2)+\frac{2\Tr G^\al(z_1)G^\al(z_2)}{n(n+2)^2}\notag
 \\
 &+ \frac{2n}{(n+2)^2}\Big(2g_n^{(1)\al}(z_1,z_2)+n^{-1}g_n^{(2)\al}(z_1,z_2)\Big)\Big]=O(n^{-1}).
\end{align}
To get $O(n^{-2})$ on the r.h.s, we need to show that the leading terms of the asymptotic expansions of the terms on the r.h.s of (\ref{nCovAa}) cancel. This is the main idea of the proof  and to implement it we use a bootstrap argument: first, we prove a weaker statement of (ii) -- (iii) and show that
\begin{align}\label{fOn-1}
%\V \gamma_n(z)=O(n),\,\,\, \,\,\,
f_n(z)=f(z)+O(n^{-1}) \,\,\,\text{and} \,\,\, f_n^{(1)}(z_1,z_2)=f(z_1)f(z_2)+O(n^{-1})
\end{align}
(note that we already have a weaker statement of (i)), and then repeating the argument and using (\ref{nCovAa}) and (\ref{fOn-1})
 we get (i) -- (iv).  Here we prove only the first part of (\ref{fOn-1}), the remaining steps are similar to the proof of Lemma 3.1 of \cite{DKL:23}.

Applying the resolvent identity, (\ref{G-G}), and (\ref{rozpis}) with $k=2$, we get
\begin{align}\label{zfn}
zf_n(z)&=-1+c-\frac{c}{m}\sum_{\al=1}^m\E\frac{1}{ A_\al}=-1+c-\frac{c}{m}\sum_{\al=1}^m\frac{1}{\E A_\al}-\frac{c}{m}\sum_{\al=1}^m\E\frac{A_\al^{\circ 2}}{A_\al(\E A_\al)^2}.
\end{align}
It follows from (\ref{|EA|}), and (\ref{VAB}), that
$\frac{c}{m}\sum_{\al=1}^m \E |A_\al^\circ|^2= O(n^{-1}).$
Also by (\ref{eega}) and (\ref{B/A}) we have
\begin{align*}
\frac{1}{m}\sum_\al \frac{1}{\E A_\al}=\frac{1}{m}\sum_\al \frac{1}{1+\tau_\al f_n^{\al}}=\frac{1}{m}\sum_\al   \frac{1}{1+\tau_\al f_n}+ O(n^{-2}).%+\frac{1}{m}\sum_{\al} \frac{\tau_\al(f_n-f_n^\al)}{(1+\tau_\al f_n^\al)(1+\tau_\al f_n)},
\end{align*}
%where by (\ref{B/A}) the second term of the r.h.s  is of order $O(n^{-2})$.
Note that if $\eta=\Im z>T$ then $|1+\tau_\al f|\ge 1-T/\eta>0$, and applying Condition A(ii) we get %for the first term on the r.h.s:
\begin{align*}
\frac{1}{m}\sum_\al\frac{1}{1+\tau_\al f_n}=\int \frac{\text{d}\sigma_n(\tau)}{1+\tau f_n}=\int \frac{\text{d}\sigma(\tau)}{1+\tau f}+ \int\frac{ \tau(f-f_n)\text{d}\sigma_n(\tau) }{(1+\tau f_n)(1+\tau f)} + O(n^{-2}).
\end{align*}
Hence
\begin{align*}
zf_n=-1+c-\int\frac{c\text{d}\sigma(\tau)}{1+\tau f}+\int\frac{c\tau(f_n-f)\text{d}\sigma_n(\tau)}{(1+\tau f_n)(1+\tau f)}+O(n^{-1}).
\end{align*}
This and  (\ref{rownanie}) lead to
\begin{align*}
z(f_n-f)=(f_n-f)\int\frac{c\tau \text{d}\sigma_n(\tau)}{(1+\tau f_n)(1+\tau f)}+O(n^{-1}).
\end{align*}
Choosing $\eta_0:=2T(c+1)$, we get for $\eta\ge \eta_0$
\begin{align*}
\left|z-\int\frac{c\tau\text{d}\sigma_n(\tau)}{(1+\tau f(z))(1+\tau f_n(z))}\right|&\geq \eta-\frac{cT}{(1-T/\eta)^2}\geq \frac{\eta}{2}\geq \frac{\eta_0}{2}.
%\left|\eta-cT\Big(1-\frac{T}{\eta}\Big)^{-2}\right|,
\end{align*}
This yields $|f_n-f|=O(n^{-1})$, and we get the first part of (\ref{fOn-1}). Then repeating the steps of the proof of Lemma 3.1 of \cite{DKL:23} and taking into account the dependence on $\tau_\al$, one can finish the proof of the Lemma \ref{lematwariancja}.
\end{proof}

It follows from Lemma \ref{lematwariancja} (i) and (\ref{gwiazdka2.4}) -- (\ref{gn-gna}), that
\begin{align}\label{Vgr}
\V(\gamma_n-\gamma_n^\al)=\V \Big(\frac{B_\al}{A_\al}\Big)=O(n^{-2}).
\end{align}
Also using (\ref{nCovAa}) and Lemma \ref{lematwariancja}, one can find the limiting covariance of the bilinear forms $K_n(z_1,z_2)$.
 Namely given $z_1$, $z_2 \in \C\backslash \R$, %, let $K_n(z_1,z_2)$ be defined in (\ref{Dn}).
%$$D_n(z_1,z_2)=\E_{Y_1,\ldots, Y_m}(\E_Y(G(z_1)Y,Y)(G(z_2)Y,Y)_Y^\circ)+\Cov_{Y_1,\ldots,Y_m}\{g_n(z_1),g_n(z_2)\}.$$
%Lemma \ref{lematwariancja} allows to show that
we have
\begin{align}\label{Dzz1}
K(z_1,z_2):=\lim_{ n \to \infty} n^2 K_n(z_1,z_2)=-2f(z_1)f(z_2)+2\frac{\Delta f}{\Delta z}+\int\frac{6c\tau^4 f^2(z_1)f^2(z_2)\text{d}\sigma(\tau)}{(1+\tau f(z_1))(1+\tau f(z_2))}
\end{align}
(see also {\cite[Theorem~3.1]{DKL:23}}.)
%Also the following statement can be proven similar to Lemmas 3.2 and 4.1 in \cite{DKL:23}.
 Finally, the next statement allows to estimate the error terms in our calculations. Its proof repeats, with appropriate modifications, the proof of Lemmas 3.2 and 4.1. We have:
\begin{lem}\label{l: var1}%[\cite{DKL:23}]
	We have uniformly in $z\in \C_{\eta_0}$ as $n\to\infty$:
	\begin{align*}
		(i)\,\,&\E_{\al}|(A_{\alpha })_\al^\circ|^p=O(n^{-p}),%\notag
		\\
		(ii)\,\,&\E|\gamma _{n}^{\circ }|^{p}=O(1),\quad%\label{gp}
		\\
		(iii)\,\,&\E|A_{\alpha }^\circ|^p=O(n^{-p}),\\
		(iv)\,\,&\V(\E_\al(A_\al^\circ)^2)=O(n^{-9/2}).
	\end{align*}
\end{lem}

Now we are ready to prove the main result.
\section{Proof of Theorem \ref{tw}}
To get Theorem \ref{tw}, it is enough to show that any linear combination of the real and imaginary parts of the trace of the resolvent, $a_1\Re\gamma_n+a_2\Im\gamma_n$, converges to a Gaussian random variable. % with zero mean and a covariance matrix $\Sigma$ defined in (\ref{mkow}). The proof of this theorem is no different from the proof of the convergence of the imaginary part of the trace of the resolvent. Therefore, below we will present the proof of the convergence of the imaginary part of the trace of the resolvent.
To this end, it suffices to consider the imaginary part $(a_1=0, a_2=1)$, the proof in the general case follows the same scheme. Thus, we need to show that $\Im(\gamma_n-\E\gamma_n)$ converges in distribution to a Gaussian random variable with zero mean and variance given by $$V(z)=\frac{1}{4}(2C(z,\overline{z})-C(z,z)-C(\overline{z},\overline{z})).$$
We use the scheme proposed in \cite{Pa-Sh} and \cite{Sh:11}.
According to the Levy's continuity theorem, it suffices to establish the  convergence of the corresponding characteristic functions. Given $z\in\C_{\eta_0}$, let
\begin{align}\label{zn}
Z_n(x)=\E e_n(x), \; \; e_n=e_n(x)=\exp\{ix\Im \gamma_n^\circ(z)\}.
\end{align}
We need to show that for every $x\in \R$,
$\lim\limits_{n \to \infty} Z_n(x)=\exp\{-x^2V(z)/2\}$.
 %Define for any test function $\varphi \in H_{s}, \ s>6.5$
%\begin{align}\label{phiy}
%\varphi_y=P_y *\varphi,
%\end{align}
%where $P_y$ is a Poisson kernel
%\begin{align}\label{Py}
%P_y(x)=\frac{y}{\pi(x^2+y^2)},
%\end{align}
%and ''$*$'' denotes the convolution. We have
%\begin{align}\label{gr}
%\lim\limits_{y \downarrow 0} ||\varphi- \varphi_y||_{s}=0.
%\end{align}
%Denote for the moment the characteristic function (\ref{zn}) by $Z_n(x)$, to make explicit its dependence on the test function. We have for any converging subsequence $\{Z_{n_j}[\phi]\}_{j=1}^{\infty}$
%\begin{align*}
%\lim\limits_{n_j \to \infty} Z_{n_j}[\varphi]=\lim\limits_{y \downarrow 0}\lim\limits_{n_j \to \infty} (Z_{n_j}[\varphi]-Z_{n_j}[\varphi_y]) +\lim\limits_{y \downarrow 0}\lim\limits_{n_j \to \infty}Z_{n_j}[\varphi_y].
%\end{align*}
%Since by (\ref{gr})
%\begin{align*}
%|Z_{n_j}[\varphi]-Z_{n_j}[\varphi_y]|\leq |x|\V\{\mathcal{N}_{n_j}[\varphi]-\mathcal{N}_{n_j}[\varphi_y]\}^{1/2}\leq C|x|||\varphi - \varphi_y||_{6.5+\sigma} \to 0,
%\end{align*}
%as $y \downarrow 0$, then
%\begin{align}\label{podciag}
%\lim\limits_{n_j \to \infty} Z_{n_j}[\varphi]=\lim\limits_{y \downarrow 0}\lim\limits_{n_j \to \infty} Z_{n_j}[\varphi_y].
%\end{align}
%Hence it suffices to find the limit of $Z_{n}(x):=Z_n(x)=\E\{e_{n}(x)\}$ with $e_{n}(x)=e^{ix\Im\gamma_n(z)}$, as
%$n \to \infty$. % It follows from (\ref{phiy})-(\ref{Py}) that
%\begin{equation}  \label{repr_N}
%\mathcal{N}_n[\varphi_ \eta ]=\frac{1}{\pi}\int
%\varphi(\mu)\Im\gamma_n(z)\textrm{d}\mu,\quad z=\mu+iy .
%\end{equation}
%This allows to write
We have
\begin{equation}
\frac{d}{dx}Z_{n}(x)=\frac{1}{2}(
\YY(z,x)- \YY(\overline z,x)),  \label{dZ=}
\end{equation}
where $\YY(z,x)=\E\gamma_n(z)e_{n}^{\circ}(x)$.
From the resolvent identity and (\ref{G-G}) (see also (\ref{zfn})), we get
\begin{align}\label{YT}
z\YY(z,x)&=-\sum
\limits_{\alpha=1}^{m} \E \frac{e_{n}^{\alpha\circ}}{A_{\alpha}}
-\sum\limits_{\alpha=1}^{m} \E \frac{e_{n}^\circ
-e_{n}^{\alpha\circ}}{A_{\alpha}}=:T_1^{(n)}+T_2^{(n)},
\end{align}
where %$A_\alpha=1+\tau_\al(G^\al Y_\al, Y_\al)$, 
 $e_n^\al=\exp\{ix\Im \gamma_n^{\al\circ}\}$ . Applying  (\ref{rozpis}) with $k=3$, we get
\begin{align*}
T_1^{(n)}&=\sum\limits_{\alpha=1}^{m} \frac{\E e_{n}^{\alpha\circ}A_{\alpha}^{\circ}}{(\E A_{\alpha})^2}-\sum\limits_{\alpha=1}^{m} \frac{\E e_{n}^{\alpha\circ}A_{\alpha}^{\circ 2}}{(\E A_{\alpha})^3}+\sum\limits_{\alpha=1}^{m} \frac{\E e_{n}^{\alpha\circ}A_{\alpha}^{\circ3}A_{\alpha}^{-1}}{(\E A_{\alpha})^3}=:T_{11}^{(n)}-T_{12}^{(n)}+T_{13}^{(n)}.
\end{align*}
It follows from Lemma \ref{l: var1} (iii), that
 $T_{13}^{(n)}=O(n^{-1})$.
Since $A_\al^\circ=(A_\al)_\al^\circ+\tau_\al g_n^{\al\circ}$ and $e_n^\al$ does not depend on $Y_\al$, we have
\begin{align*}
	\E e_n^{\al\circ}A_\al^{\circ 2}=\E(e_n^{\al\circ}\E_\al(A_\al)_\al^{\circ 2})+\tau_\al^2\E e_n^{\al\circ}(g_n^{\al\circ})^2.
\end{align*}
Applying Lemma \ref{l: var1} (iv) to the first term and Lemma \ref{lematwariancja} (i) to the second term we get
\begin{align*}
	|\E e_n^{\al\circ} A_\al^{\circ 2}|\lesssim (\V(\E_\al(A_\al)_\al^{\circ 2}))^{1/2} + \V g_n^\al= O(n^{-9/4}),
\end{align*}
so that $T_{12}^{(n)}=O(n^{-1/4})$. Consider now $T_{11}^{(n)}$. Since $e_{n}^{\alpha}$ does not depend on $Y_{\alpha}$,
\begin{align*}
\E e_{n}^{\alpha\circ}A_{\alpha}=\E e_{n}^{\alpha\circ}\E_{\alpha}A_{\alpha}=\tau_\alpha n^{-2}\E e_{n}^{\alpha\circ}\gamma_{n}^{\alpha}=\tau_\alpha n^{-2}\YY(z,x)+R_n,
\end{align*}
where $$R_n=\tau_\alpha n^{-2}\E(e_{n}^{\alpha\circ} \gamma_n^{\alpha}-e_{n}^{\circ}\gamma_n)=\tau_\al n^{-2}\E((e_{n}^{\alpha}-e_{n}) \gamma_n^\circ+e_{n}^{\alpha}(\gamma_n^\al-\gamma_n)^\circ).$$
Using the Taylor series expansion, we get
\begin{align}\label{eyn-eyna}
	e_{n}-e_{n}^{\alpha}=ixe_{n}^{\alpha}\Im(\gamma_n-\gamma_n^{\alpha})^{\circ}+O(|\Im(\gamma_n-\gamma_n^{\alpha})^{\circ}|^2).
\end{align}
This and (\ref{Vgr}) leads to
$|R_n|=O(n^{-3})$.
Hence $$\E e_n^{\al\circ}A_\al=\tau_\al n^{-2} \YY(z,x)+ O(n^{-3})$$ and
\begin{align*}
T_{11}^{(n)}=\sum\limits_{\alpha=1}^{m} \frac{\E e_{n}^{\alpha\circ}A_{\alpha}^{\circ}}{(\E A_{\alpha})^2}=\YY(z,x)\int\frac{c\tau\text{d}\sigma(\tau)}{(1+\tau f(z))^2}+O(n^{-1}).
\end{align*}
Summarising we get
\begin{align}\label{T1nw}
T_{1}^{(n)}=\YY(z,x)\int\frac{c\tau\text{d}\sigma(\tau)}{(1+\tau f(z))^2}+ O(n^{-1/4}).
\end{align}
Consider $T_2^{(n)}$ of (\ref{YT}). It follows from (\ref{eyn-eyna})
%\begin{align*}
%e_{n}-e_{n}^{\alpha}&=\frac{ixe_{n}^{\alpha}}{\pi}\int\varphi(\lambda_1)\Im(\gamma_n-\gamma_n^{\alpha})^{\circ}(z)\textrm{d}\lambda_1\nonumber\\&+O\Big(\Big|\int\frac{1}{\pi}\varphi(\lambda_1)\Im(\gamma_n-\gamma_n^{\alpha})^{\circ}(z)\textrm{d}\lambda_1\Big|^2\Big),
%\end{align*}
and (\ref{gn-gna}), that
\begin{align*}
\E A_{\alpha}^{-1}(z)(e_{n}-e_{n}^{\alpha})^{\circ}&=ix\E e_{n}^{\alpha}(A_{\alpha}^{-1}(z))^{\circ}\Im(B_{\alpha}A_{\alpha}^{-1})^{\circ}(z)+O(R_n^{(1)}),
\end{align*}
where applying Cauchy-Schwarz inequality, (\ref{gwiazdka2.4}) and Lemma \ref{l: var1} (iii) one can get
\begin{align*}
|R_n^{(1)}|&\le \E|(A_\al^{-1})^\circ||(\gamma_n-\gamma_n^\al)^\circ|^2\le\Big(\E|(A_\al^{-1})^\circ|^2\E \big|\big(B_\al/A_\al\big)^\circ\big|^4\Big)^{1/2}\\\notag&\lesssim\big ( \V A_\al)^{1/2} (\E|A_\al^\circ|^4+\E|B_\al^\circ|^4\big)^{1/2}=O(n^{-3}).
\end{align*}
Hence
\begin{align*}
\E A_{\alpha}^{-1}(z)(e_{n}-e_{n}^{\alpha})^{\circ}&=\frac{x}{2}\E e_{n}^{\alpha}(A_{\alpha}^{-1}(z))^{\circ}[(B_{\alpha}A_{\alpha}^{-1})^{\circ}(z)-(B_{\alpha}A_{\alpha}^{-1})^{\circ}(\overline{z})]+O(n^{-3}).
\end{align*}
Applying twice (\ref{rozpis}) with $k=2$, and then using Lemma  \ref{l: var1} (iii), we get for $z_1,z_2\in \C_{\eta_0}$
\begin{align*}
&\E e_{n}^{\alpha}(A_{\alpha}^{-1}(z_1))^{\circ}(B_{\alpha}A_{\alpha}^{-1})^{\circ}(z_2)\notag\\&=\E e_n^\al\Big(-\frac{A_\al^\circ}{(\E A_\al)^2}+\frac{A_\al^{\circ 2}}{A_\al(\E A_\al)^2}\Big)(z_1)\Big(\frac{B_\al^\circ}{\E A_\al}-\frac{(B_\al A_\al^\circ)^\circ}{(\E A_\al)^2}+\frac{(B_\al A_\al^{\circ 2})^\circ}{A_\al(\E A_\al)^2}\Big)(z_2)+O(n^{-3})\notag\\&
=-\frac{\E e_{n}^{\alpha}A_{\alpha}^{\circ}(z_1)B_{\alpha}^{\circ}(z_2)}{(\E A_{\alpha}(z_1))^2\E A_{\alpha}(z_2)}+\frac{\E B_\al(z_2) \E e_{n}^{\alpha}A_{\alpha}^{\circ}(z_1)A_{\alpha}^{\circ}(z_2)}{(\E A_{\alpha}(z_1))^2(\E A_{\alpha}(z_2))^2}+O(n^{-3}).
\end{align*}
It follows from Lemma \ref{l: var1} (iv), that
\begin{align*}
\E e_{n}^{\alpha}A_{\alpha}^{\circ}(z_1)A_{\alpha}^{\circ}(z_2)=\E e_{n}^{\alpha}\E A_{\alpha}^{\circ}(z_1)A_{\alpha}^{\circ}(z_2)+O(n^{-9/4}).
\end{align*}
Also, by (\ref{Vgr}) and (\ref{eyn-eyna}), we have
$$|\E e_n^\al-Z_n|=|\E(e_n^\al-e_n)|=O(n^{-1}),$$
where $Z_n$ is defined in (\ref{zn}). Hence, taking into account that $B_{\alpha}=\frac{\partial}{\partial z}A_{\alpha}$, we get  for $z_1\neq z_2$
\begin{align*}
\E e_{n}^{\alpha}(A_{\alpha}^{-1}(z_1))^{\circ}(B_{\alpha}A_{\alpha}^{-1})^{\circ}(z_2)&=\frac{-Z_{n}(x)}{(\E A_{\alpha}(z_1))^2}\frac{\partial}{\partial z_2}\frac{\E A_{\alpha}^{\circ}(z_1)A_{\alpha}^{\circ}(z_2)}{\E A_{\alpha}(z_2)}+O(n^{-3})\nonumber\\&=\frac{-Z_{n}(x)}{(1+\tau_\al f(z_1))^2}\frac{\partial}{\partial z_2}\frac{\tau_\al^2K(z_1,z_2)}{1+\tau_\al f(z_2)}+O(n^{-3}),
\end{align*}
 where $K(z_1,z_2)$ is defined in (\ref{Dzz1}). % $\Delta f=f(z)-f(z_1)$ and $\Delta z=z-z_1$.
Hence
\begin{align*}
T_2^{(n)}=-\frac{cxZ_{n}(x)}{2}\int\frac{\tau^2}{(1+\tau f(z))^2}\lim\limits_{z_2\to z}\frac{\partial}{\partial z_2}\bigg[\frac{K(z,z_2)}{1+\tau f(z_2)}-\frac{K(z,\overline{z}_2)}{1+\tau f(\overline{z}_2)}\bigg]\text{d}\sigma(\tau)+O(n^{-1}).
\end{align*}
This,  (\ref{YT}), and (\ref{T1nw}) yield
\begin{align*}
Y_n(z,x)&=\bigg(\int\frac{c\tau\text{d}\sigma(\tau)}{(1+\tau f(z))^2}-z\bigg)^{-1}\\
&\times \frac{xZ_{n}(x)}{2}\int\frac{c\tau^2}{(1+\tau f(z))^2}\lim\limits_{z_2\to z}\frac{\partial}{\partial z_2}\bigg[\frac{K(z,z_2)}{1+\tau f(z_2)}-\frac{K(z,\overline{z}_2)}{1+\tau f(\overline{z}_2)}\bigg]\text{d}\sigma(\tau)+O(n^{-1/4}).
\end{align*}
Equation (\ref{rownanie}) and the resulting equalities
\begin{align*}
\bigg(\int\frac{c\tau\text{d}\sigma(\tau)}{(1+\tau f)^2}-z\bigg)^{-1}=f^\prime/f, \quad
\frac{\Delta z}{\Delta f}=\frac{1}{f(z_1)f(z_2)}-\int\frac{c\tau^2\text{d}\sigma(\tau)}{(1+\tau f(z_1))(1+\tau f(z_2))}
\end{align*}
 allow us to show that
$$\bigg(\int\frac{c\tau\text{d}\sigma(\tau)}{(1+\tau f(z))^2}-z\bigg)^{-1} \int\frac{c\tau^2}{(1+\tau f(z))^2}\frac{\partial}{\partial z_2}\frac{K(z,z_2)}{1+\tau f(z_2)}\text{d}\sigma(\tau)=C(z,z_2),$$
where $C(z,z_2)$ is defined in (\ref{CZZ}) (see also \cite{GLPP:13}.) Hence
\begin{align*}
Y_n(z,x)=\lim\limits_{z_2\to z} \frac{xZ_{n}(x)}{2}\big[C(z,z_2)-C(z,\overline{z}_2)\big]+O(n^{-1/4}).
\end{align*}
 This and (\ref{dZ=}) lead to
\begin{align*}
\frac{\partial}{\partial x}Z_{n}(x)=-xV(z)Z_{n}(x)+O(n^{-1/4}),
\end{align*}
 and finally we get
\begin{align*}
\lim_{n\to\infty}Z_{n}(x)=\exp\{-x^2V(z)/2\}.
\end{align*}
%Now we take into account (\ref{podciag}), allowing us to pass to the limit $y\downarrow 0$, and obtain (\ref{t:clt}).
The theorem is proved.\\

\textbf{Acknowledgments.} The author would like to thank Anna Lytova for introducing to the problem and for the fruitful discussions.
The project was supported by grant nr 2018/31/B/ST1/03937 National Science Centre, Poland.

\address
\end{document}